\documentclass[a4paper]{amsart}

\title[triple-loop networks with arbitrarily many MDD's]
{Triple-loop networks with arbitrarily many \\ minimum distance diagrams}
\thanks{Research partially supported by  the Spanish Ministry of Education and Science, grant number
 MTM2005-08618-C02-02.}

\author{Pilar Sabariego
        \and
        Francisco Santos}
        \address{Departamento de Matem\'aticas, Estad\'{\i}stica y
        Computaci\'on, Universidad de Cantabria, Santander, Spain}
        \email{pilar.sabariego@unican.es, francisco.santos@unican.es.}

\usepackage{latexsym}
\usepackage{amssymb}
\usepackage{amsmath}
\usepackage[dvips]{graphicx}

\usepackage{color}

\usepackage{multirow}

\newtheorem{theorem}{Theorem}[section]

\newtheorem{corollary}[theorem]{Corollary}
\newtheorem{lemma}[theorem]{Lemma}
\newtheorem{definition}[theorem]{Definition}

\newtheorem{remark}[theorem]{Remark}
\newtheorem{example}[theorem]{Example}

\newcommand\Z{\mathbb Z}
\newcommand\R{\mathbb R}
\newcommand\N{\mathbb N}

\newcommand{\pos}{\operatorname{pos}}

\renewcommand{\mod}{\operatorname{mod}\ }

\begin{document}

\begin{abstract}
\emph{Minimum distance diagrams} are a way to encode the diameter and routing information of multi-loop networks. For the widely studied case of double-loop networks, it is known that each network has at most two such diagrams and that they have a very definite form (``$L$-shape'').

In  contrast, in this paper we show that there are triple-loop networks with an arbitrarily big number of associated minimum distance diagrams. For doing this, we build-up on the relations between minimum distance diagrams and monomial ideals.
\end{abstract}

\maketitle

\section{Introduction}

Multi-loop networks have been widely used in the computer and network architecture literature, as a simple, yet efficient, way of organizing multi-module memory services. Their mathematical study was initiated in~\cite{Wong-Coppersmith-1974}, where the problem of finding the network parameters that minimize the diameter (and/or the average distance) for networks of given size and degree was posed.

\begin{definition}
A \emph{multi-loop network} of size $N$ and steps $s_{1}, \dots ,s_{r}$ is a directed graph with nodes $V=\{0,1,\dots,N-1\} = \Z_{N}$ and an arc $i\to i+s_l$ for every $i\in V$ and every $s_l$, i.e.,
\[
i \rightarrow i+s_{l} \mod  N.
\]
We denote this network by $C_{N}(s_{1}, \dots ,s_{r})$. 
\end{definition}

In other words, $C_{N}(s_{1}, \dots, s_{r})$ is a Cayley digraph of the ciclic group $\Z_{N}$ with respect to $\{s_{1}, \dots, s_{r}\}$.

One convenient way to encode the routing information in these networks is assigning to each vertex $i\in\Z_N$ an integer non-negative vector $\mathbf a=(a_1,\dots, a_r)$ such that 
\[
a_1s_1 + \cdots + a_rs_r = i \mod N.
\]
Then, one can go from node 0 to node $i$ by traversing $a_i$ nodes of length $s_i$, for each $i$. The order is irrelevant. Also, since the network is vertex-transitive. the same is valid for any pair of vertices $j$ and $j+i (\mod N)$.

A \emph{minimum distance diagram} (MDD for short) for the network $C_N(s_1,\dots, s_r)$ is just this information, except the assumption is made that the path taken for each node $i$ has minimal length. 
In particular, from a minimum distance diagram we can calculate the diameter and the average distance of the circulant digraph.

\begin{figure}[htb]
      \includegraphics[width=6cm]{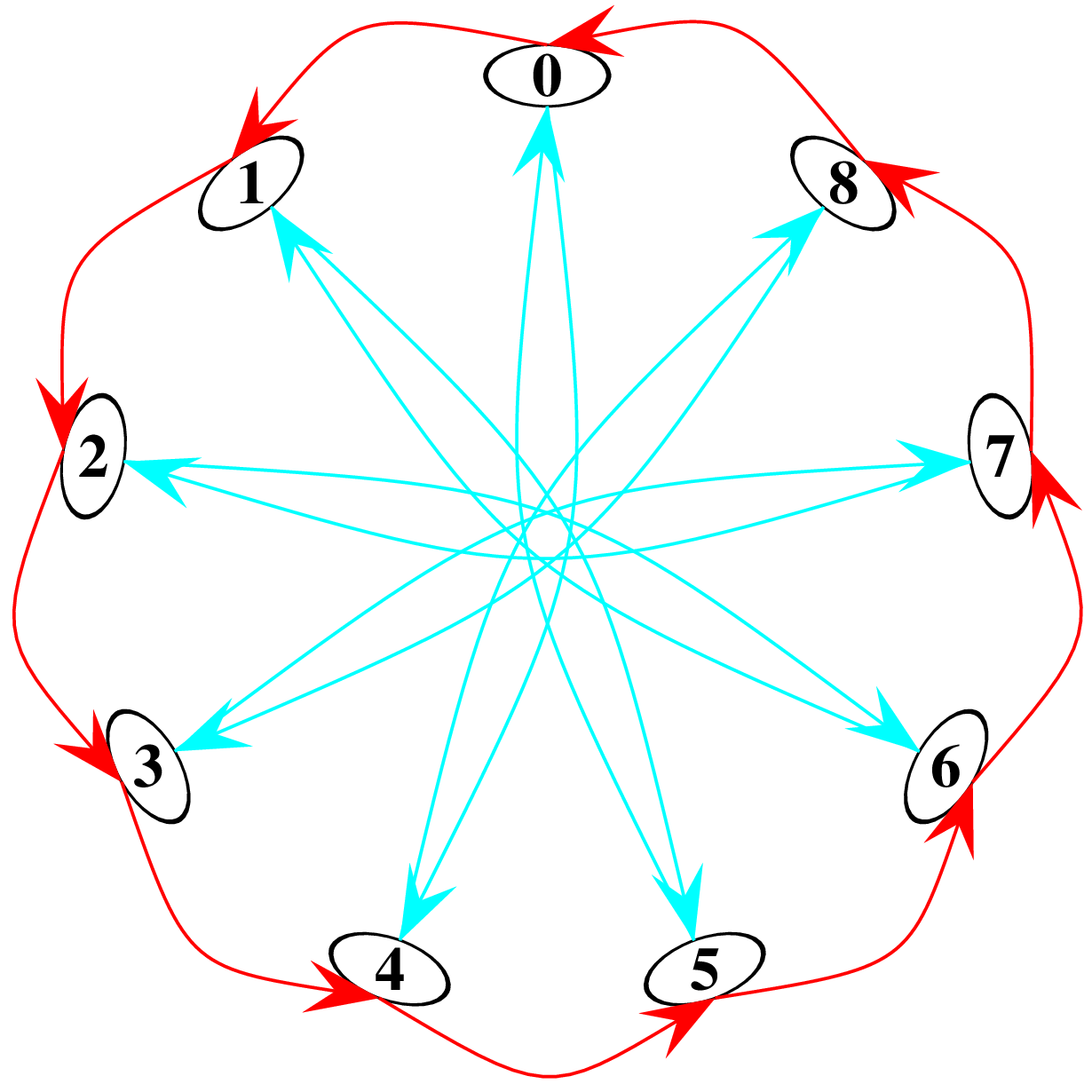}
\qquad 
\begin{minipage}[b]{5.5cm} 
      \includegraphics[width=2.4cm]{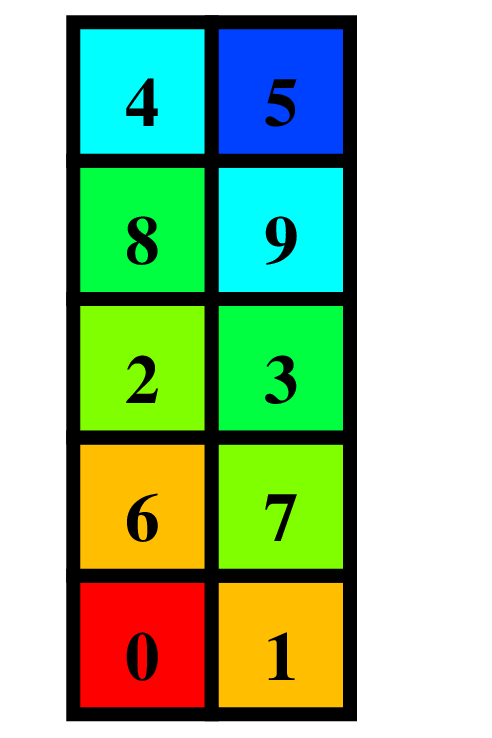}
\\
      \includegraphics[width=4.9cm]{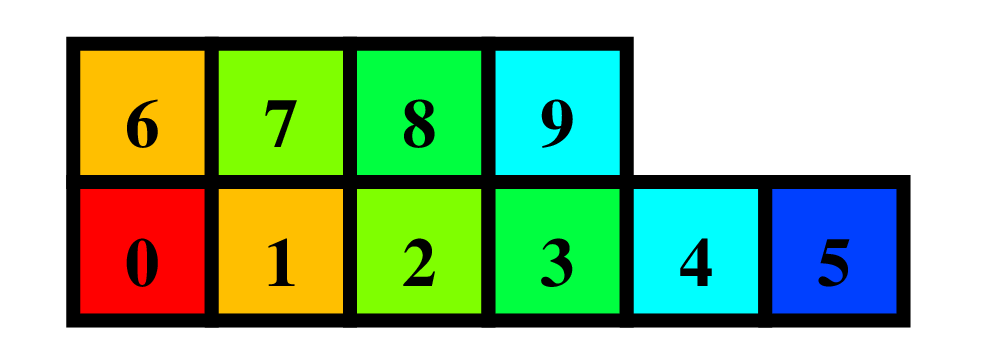}
\vspace{0.2cm}
\end{minipage} 
\caption{The network $C_{9}(1,4)$ and its two MDD's ($L$-shapes).}
\label{fig:digrafoyMDDs}
\end{figure}

An example is in Figure~\ref{fig:digrafoyMDDs},  for the network $C_{9}(1,4)$. The left part of the picture represents the network itself and the right are two minimum distance diagrams of it, in their customary graphical representation. In the first diagram, we have chosen as minimum path to vertex $7$ the one that takes three steps of type $s_1=1$ (horizontal steps in the diagram) and one step of type $s_2=4$ (vertical step). But we can as well choose four steps of type $s_2$ since $1+1+1+4=4+4+4+4=7 (\mod 9)$.

As seen in the picture, the minimum distance diagrams are multidimensional ``stacks of cubes''. They tesselate the space by the action of the following natural lattice associated to the  network.
\begin{equation}
\mathcal{L}:=\{(a_{1}, \dots, a_{r}) \in \Z^{r}: a_{1}s_{1}+\cdots+a_{r}s_{r}\equiv 0 \mod N\}.
\label{eq:lattice}
\end{equation}

Minimum distance diagrams appear frequently in the literature on multi-loop networks, although normally without a \emph{definition} of what an MDD is in general. Rather, since the interest is in solving the routing problem, the authors consider one particular diagram, and concentrate on the algorithm to obtain it, or on studying its shape, etc. 
Most authors~\cite{Aguilo-Fiol-1995, Aguilo-etal-1997, Aguilo-Miralles-2004, Chen-Hwang-1988, Chen-etal-2003, Du-etal-1985, Fiol-etal-1987, Hwang-2001, Hwang-2003} always ``break ties'' lexicographically whenever there are two minimal paths to a vertex, so that every multi-loop network has a unique MDD for each prescribed ordering of the parameters $s_1,\dots,s_r$.

In our definition (see Section 2), besides asking each individual path to have minimum length among those to a given vertex, we include a technical condition that is implicitly present in all previous work and which can be rephrased saying that the diagram is the complement of an ideal in $\N^r$.

It is known~\cite{Wong-Coppersmith-1974,Hwang-2001} that the MDD's of double-loop networks have a very precise form for which they are called \emph{L-shapes} (see Figure~\ref{fig:digrafoyMDDs}). Aguil\'o and Miralles~\cite{Aguilo-Miralles-2004} have shown that for each double-loop network there are at most two such $L$-shapes that are MDD's for it. We give a new proof of this in Lemma~\ref{lemma:2MDDs}. From this characterization of the shape of MDD's it is easily derived that a double-loop network with diameter $D$ cannot have more than $N\simeq D^2/3 +O(D)$ nodes. Networks that achieve this bound are known.

\begin{figure}[htb]
\centering
 \includegraphics[width=5.7cm]{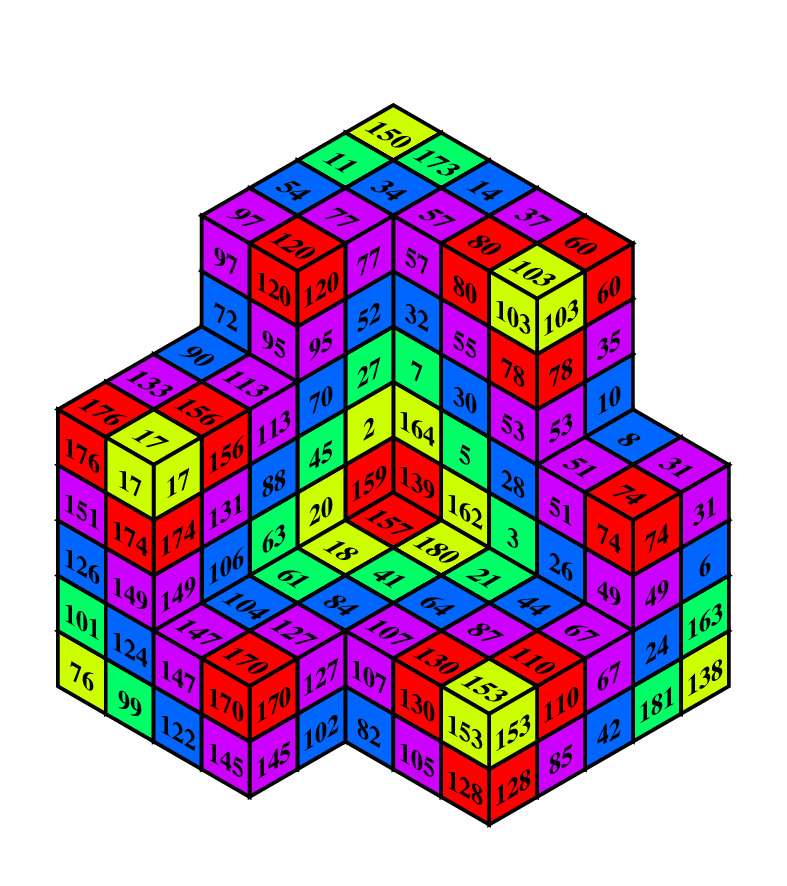}
 \caption{The ``hyper-L with parameters $l=7$, $m=3$ and $n=2$''~\cite{Chen-etal-2003}. It is a minimum distance diagram for the circulant digraph $C_{182}(43,23,25)$. 
 }
\label{fig:hyper-Ls}
\end{figure}

In order to construct triple-loop networks with low diameter, Aguil\'o et al.~\cite{Aguilo-etal-1997,Aguilo-1999} have considered similar nicely shaped MDD's for them, 
the so-called  \emph{hyper-L tiles} (see Figure~\ref{fig:hyper-Ls}). These are  MDD's of certain triple-loop networks $C_N(s_1,s_2,s_3)$ with diameter $D$ satisfying\footnote{
Observe that triple-loop networks with $N\simeq D^3/27$ are trivial to construct, and that, by a simple volume argument on its MDD, every triple loop
 network has $N\le \binom{D+3}{3}\simeq D^3/6$. A better upper 
 bound of $N\le (D+3)^3/(14-3\sqrt{3})\simeq 0.11 D^3$ was given by Hsu and Jia~\cite{Hsu-Jia-1994}.}
\[
N\ge \frac{2}{27} D^3 + O(D^2).
\]
But it was shown in \cite{Chen-etal-2003, Chen-etal-2006} that these hyper-L MDD's exist only for very special parameters $N$, $s_1, s_2$ and $s_3$ of the network.

Our initial goal in this work was to get an upper bound for the number of MDD's of triple-loop networks. But  the truth is that a global bound does not exist. This indicates it is certainly a difficult task to characterize them:

\begin{theorem}[Theorem~\ref{thm:many}]
\label{thm:many-intro}
If $q\in \N$ is not a multiple of three and we let $k=2+q+q^2$, the triple-loop network $C_{k(k-1)}(1+k,1+kq,1+kq^2)$ has exactly $3(q+2)$ minimum distance diagrams.
\end{theorem}

\begin{figure}[htb]
\vspace{1cm}
\begin{center}
\begin{tabular}{ccc}

\begin{minipage}{2cm}
\begin{center}
      \includegraphics[width=2cm]{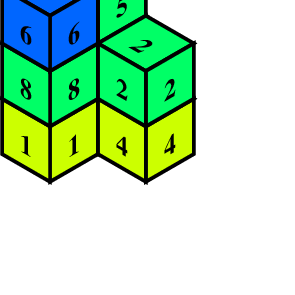}
\end{center}
\end{minipage}  
&
\begin{minipage}{2cm} 
\begin{center}
      \includegraphics[width=2cm]{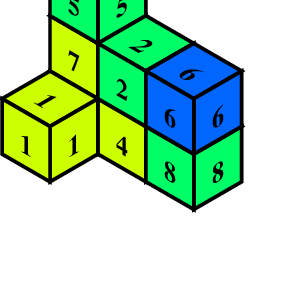}
\end{center}
\end{minipage} 
&
\begin{minipage}{2cm} 
\begin{center}
      \includegraphics[width=2cm]{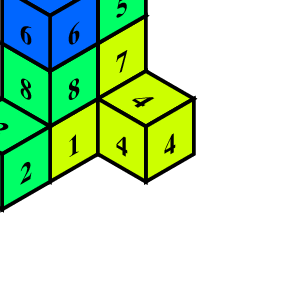}
\end{center}
\end{minipage} 
\\
\begin{minipage}{2cm} 
\begin{center}
      \includegraphics[width=2cm]{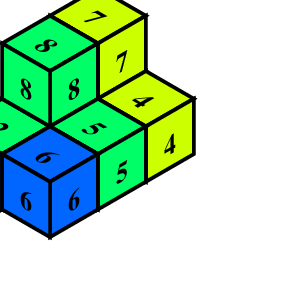}
\end{center}
\end{minipage} 
&
\begin{minipage}{2cm} 
\begin{center}
      \includegraphics[width=2cm]{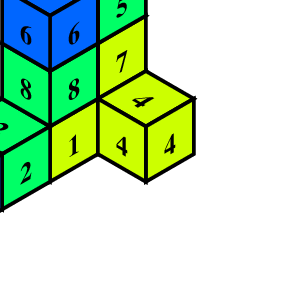}
\end{center}
\end{minipage} 
&
\begin{minipage}{2cm} 
\begin{center}
      \includegraphics[width=2cm]{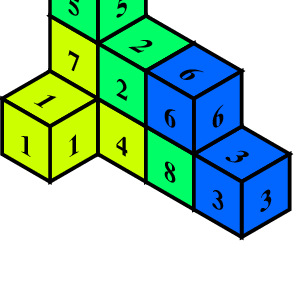}
\end{center}
\end{minipage} 
\\
\begin{minipage}{2.5cm} 
\begin{center}
      \includegraphics[width=2cm]{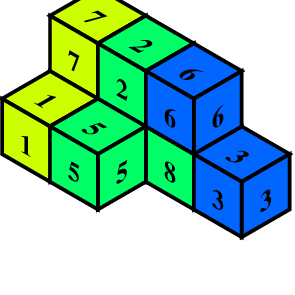}
\end{center}
\end{minipage} 
&
\begin{minipage}{2.5cm} 
\begin{center}
      \includegraphics[width=2cm]{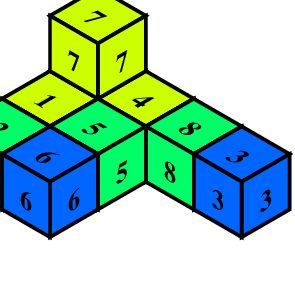}
\end{center}
\end{minipage} 
&
\begin{minipage}{2.5cm} 
\begin{center}
      \includegraphics[width=2cm]{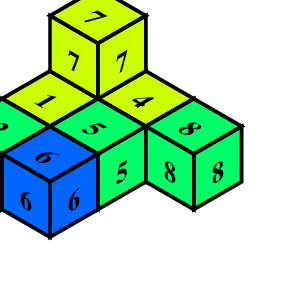}
\end{center}
\end{minipage} 
\\
\end{tabular}
\end{center}
\caption{Circulant Digraph $C_{9}(1,4,7)$ and its 9 associated MDD's.}
\label{fig:3Pasosy9MDDs}
\end{figure}

For a small example consider the following:  $C_9(1,4,7)$, has nine different MDD's, shown in Figure~\ref{fig:3Pasosy9MDDs}.

Our proof of Theorem~\ref{thm:many-intro}
is based in the interpretation of minimum distance diagrams in terms of initial ideals of a certain lattice ideal. More precisely: 
\begin{theorem}[G\'omez et al., 2006]
Let $\mathcal L$ denote the lattice of equation (\ref{eq:lattice}). Then,
the complement of any monomial graded initial ideal of the lattice ideal of $\mathcal L$ is an MDD for
the network $C_N(s_1,\dots,s_r)$.
\end{theorem}

See Section 3 for more details on lattice ideals and the role they play in multi-loop networks. Following the terminology used in the theory of toric ideals~\cite{Sturmfels-1996} we call the MDD's that can be obtained as initial ideals~\emph{coherent}. We show small examples of non-coherent MDD's for quadruple loop networks (Example~\ref{ex:non-coherent}) but do not know whether they exist for triple-loop networks. In particular, all the MDD's obtained in Theorem~\ref{thm:many-intro} are coherent.

The interpretation of coherent MDD's as initial ideals relates our result to the following statement from~\cite{Sturmfels-etal-1995}: ``There are lattice ideals in dimension three with arbitrarily large Gr\"obner bases''. The size of a Gr\"obner basis for the lattice ideal of a multi-loop network is related to the ``combinatorial complexity'' of the associated MDD. For example, the L-shape property of MDD's in double-loop networks is a consequence of the following result, also from~\cite{Sturmfels-etal-1995}: ``Gr\"obner bases for two-dimensional lattice ideals have at most three elements''.


 The rest of the paper is divided into four sections: In Section~\ref{sec:digraph-mdds} we give our precise definition of minimum distance diagrams. In Section~\ref{sec:mdds-ideals} we recall and extend the above mentioned result from~\cite{Gutierrez-etal-2006} that relates minimum distance diagrams to initial ideals of a lattice ideal (Theorem~\ref{thm:elimination-ideals}).
In Section~\ref{sec:ManyMDD's}, we concentrate on triple-loop networks and show how their number of MDD's can be bounded above by the cardinalities of the Hilbert bases of a two-dimensional \emph{homogeneous lattice} associated to the network. Finally, in Section~5 we show that this bound is tight in some cases, and use this to construct triple-loop networks with arbitrarily many associated minimum distance diagrams.

\section{Multi-loop networks and MDD's}
\label{sec:digraph-mdds}

The routing problem in a multi-loop network $C_N(s_1,\dots,s_r)$ (that is, finding the minimum path between two given vertices $i$ and $j$) can be rephrased as the following diophantine programming problem:
minimize $|a_1+\cdots + a_r|$ such that 
$j-i= a_1s_1+\cdots a_rs_r \, (\mod N)$ and $(a_1,\dots,a_r)\in \N^r$ (where, by convention, we take $\N=\{0,1,2,\dots\}$). Indeed, a path from $i$ to $j$ will always consist of a certain number $a_i\in \N$ of arcs of each type $s_i$, and the relative order in which steps are done does not affect the length of the path.

There is no loss of generality in assuming that $i=0$ (multi-loop networks are vertex-transitive), and a minimum path can be represented simply by the vector $(a_1,\dots,a_r)\in \N^r$. As said in the introduction, this suggest the following definition:

\begin{definition}
\label{def:MDDs}
A \emph{minimum distance diagram} (MDD for short) for a multi-loop network $C_{N}(s_{1}, \dots, s_{r})$ is any map $D: \Z_{N} \to \N^r$
such that:
\begin{enumerate}
\item For every $i\in  \Z_{N}$, $D(i) = (a_1,\dots,a_r)$ satisfies
$i=a_1s_1+ \cdots + a_rs_r\, (\mod N)$
and  $\|D(i)\|_1$ is minimum among all the vectors in $\N^r$ with that property.

\item For every $i$ and for every vector $b\in \N^r$ that is coordinate-wise smaller than $D(i)$ we have $b=D(j)$ for some $j$. (Of course, for this to be possible we must have  $j=b_1 s_1 + \cdots + b_rs_r$).
\end{enumerate}
\end{definition}

Property (1) in the definition says just that the map $D$ gives a solution to the routing problem 
for each  $i$. The second condition is a ``compatibility'' or ``consistency'' condition on the solutions for different values of $i$. It states that 
if one of the paths from vertex $i$ to vertex $j$ specified by the MDD passes through a vertex $k$, then the two subpaths from $i$ to $k$ and from $k$ to $j$ are also among those specified in the MDD. 
This condition is not required by, for example, G\'omez et al.~\cite{Gutierrez-etal-2006}, but holds for all the MDD's considered in the literature and is sometimes implicitly assumed.

\begin{remark}
\rm
It is clear that knowing  the image $D( \Z_{N})$ of $D$ is enough to describe $D$. 
For this reason
we will often abuse language and call minimum distance diagram the image of $D$. This is also done in the literature, where an MDD is usually characterized by ``its shape'', and it justifies the name ``diagram'' for it.

In this sense, MDD's admit (for small $r$) a nice graphical representation as a ``stack of labeled boxes'': boxes represent elements of $ \mathbb{N}^{r}$ and they are labeled by the numbers $0,\dots, N-1$. 
\end{remark}

\begin{figure}[htb]
\begin{minipage}{6cm}
      \includegraphics[width=5.5cm]{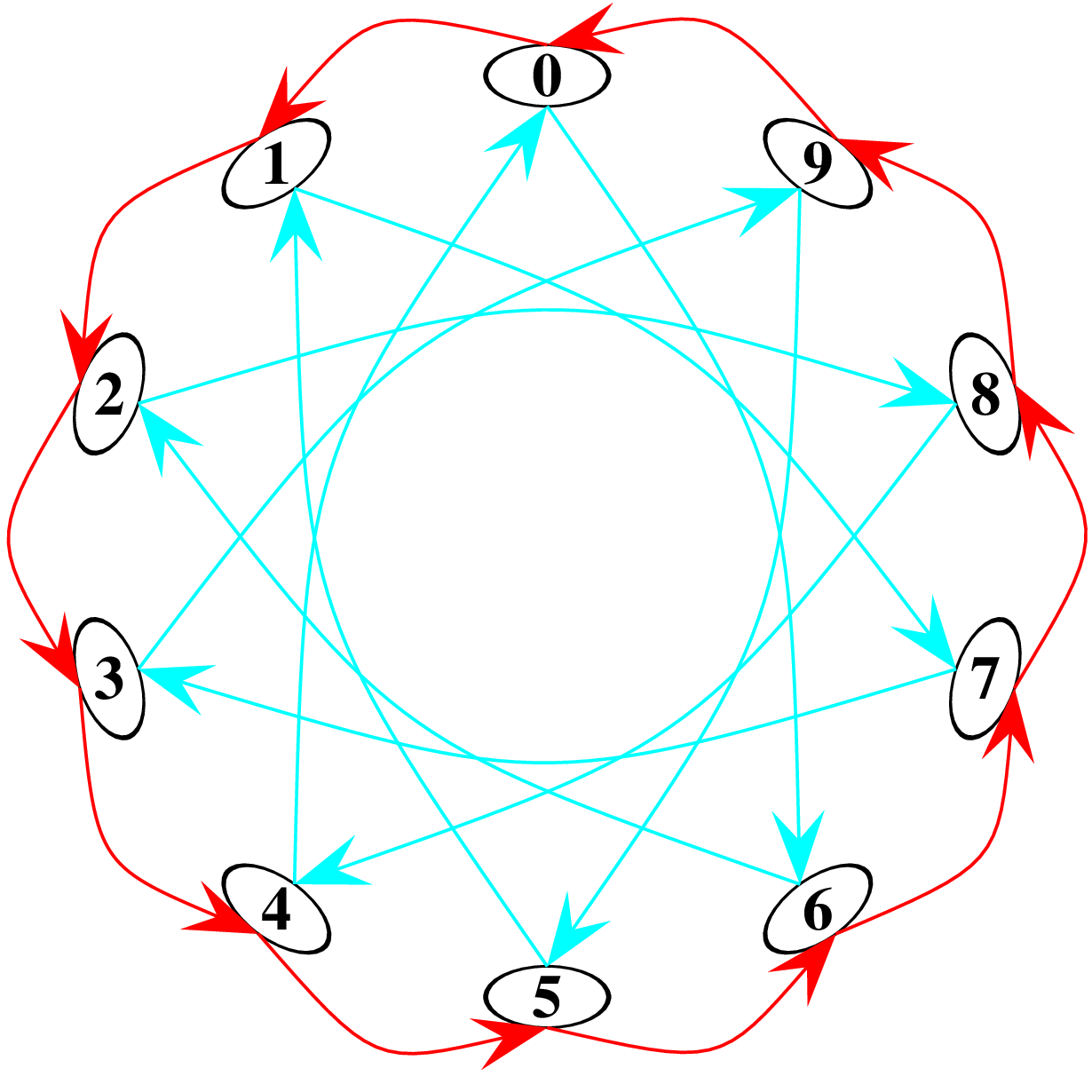}
\end{minipage}  
\begin{minipage}{5cm} 
      \includegraphics[width=4.5cm]{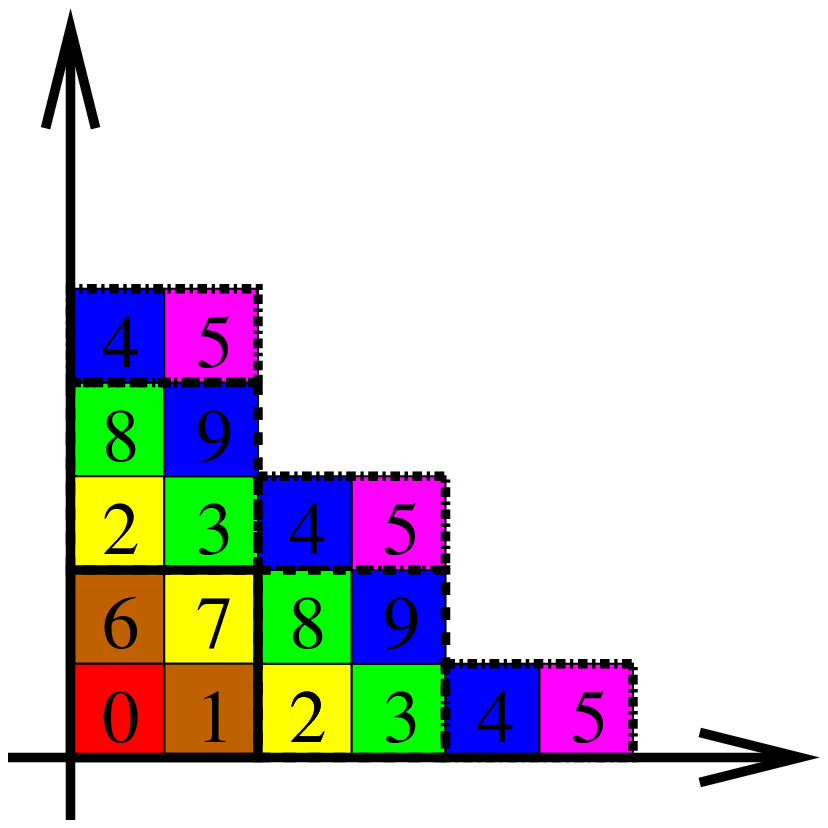}
\end{minipage} 
\caption{The network $C_{10}(1,6)$ and part of the tesselation of the plane given by the routing map associated to it.}
\label{fig:digrafoyRoutingMap}
\end{figure}

\begin{example}
\rm
To understand  better what the second condition means, let us consider the double loop network $C_{10}(1,6)$, drawn in Figure~\ref{fig:digrafoyRoutingMap}. The right picture of the figure  shows part of its routing map. In fact, it shows, for each node $c\in \Z_{N}$, \emph{all} the minimum routes from $0$ to $c$ in the network. There is a unique one for $c \in \{0,1,6,7\}$, there are two for $c \in \{2,3,8,9\}$ and there are three for $c \in \{4,5\}$. That is, this network has exactly $2^{4}3^{2}=144$ ``diagrams" that verify the first condition of an MDD.

But most of these diagrams are not very natural. If we choose to go from $0$ to $2$ by two steps of length $s_1=1$, instead of two steps of length $s_2=6$, it seems natural to go from $0$ to $8$ by two steps of length $1$ and one of length  $6$ rather than by three steps of length $6$. This is what the second condition asks for.
\end{example}

\begin{figure}[htb]
\begin{minipage}[b]{1.5cm}
\includegraphics[height=3cm]{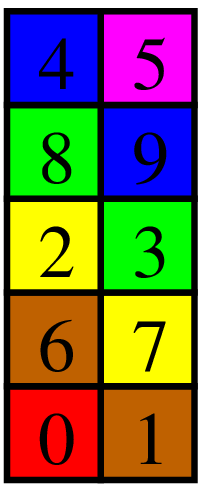}
\end{minipage}
\qquad
\begin{minipage}[b]{4cm}
\includegraphics[width=3.6cm]{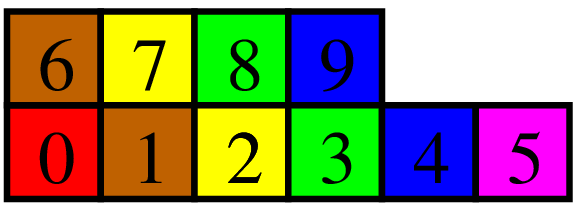}
\end{minipage}
 \caption{Two MDD's for the network $C_{10}(1,6)$.}
 \label{fig:2L_shapes}
\end{figure}

The case $r = 2$ (``double loop networks'') is the most studied type of multi-loop networks and there are several characterizations and studies of them and their MDD's, \cite{Aguilo-Fiol-1995,Chen-Hwang-1988, Du-etal-1985, Fiol-etal-1987}. 
In particular, it is known that MDD's have a very nice shape usually called  \emph{L-shape} and that, moreover, every double-loop network has at most two such $L$-shapes. For example, going back to our example of $C_{10}(1,6)$, in Figure~\ref{fig:2L_shapes} we see the only two MDD's. More precisely, the choice for $D(2)$, which can be equal to either $(2,0)$ or $(0,2)$, fixes the rest of the MDD.

\begin{lemma}[Aguil\'o and Mirall\'es, 2004]
\label{lemma:2MDDs}
Every double-loop network has exactly two MDD's.
\end{lemma}

\begin{proof}
Let $C_N(s_1,s_2)$ be our network. If the network admits more than one MDD then, in particular, there must be some $i\in \Z_N$ such that there is a choice for $D(i)$. That is, there are two vectors ${\mathbf a}=(a_1,a_2)$ and ${\mathbf b}=(b_1,b_2)$ in $\N^2$ with the same $L_1$-norm and with $a_1s_1+a_2s_2=b_1s_1+b_2s_2$. 

It is easy to see, also, that if we choose ${\mathbf a}$ and ${\mathbf b}$ in that conditions \emph{and with minimum possible $L_1$-norm} $||{\mathbf a}|| = ||{\mathbf b}|| = l$, then they must be equal to $(l,0)$ and $(0,l)$, respectively. Indeed, if both $a_1$ and $b_1$ are positive, we can subtract $(1,0)$ from ${\mathbf a}$ and ${\mathbf b}$ and if both $a_2$ and $b_2$ are positive, we can subtract $(0,1)$. We claim that which of $(l,0)$ and $(0,l)$ belongs to a particular MDD completely determines the rest of the MDD: if for some other $D(j)$ we have two (or more) choices, say $(c_1, c_2)$ and $(d_1,d_2)$, then 
$(c_1, c_2) - (d_1,d_2)$ is an integer multiple of $(l,-l)$, so one of them is incompatible with $(l,0)$ and the other with $(0,l)$.
\end{proof}

 \section{MDD's and monomial ideals}
 \label{sec:mdds-ideals}
 In this section we recall a result of~\cite{Gutierrez-etal-2006} relating MDD's of a multi-loop network with initial ideals of a certain lattice ideal. We also offer an algorithm to compute this ideal, different from the one in~\cite{Gutierrez-etal-2006}.

 Let $\mathbb{K}$ be an arbitrary field and let $\mathbb{K}[x_{1},\dots,x_{r}]$ be the polynomial ring in the variables $x_{1},\dots, x_{r}$. As customary, monomials of $\mathbb{K}[x_{1},\dots,x_{r}]$ are identified with vectors of $\N^{r}$ in the following natural  way: 
 \begin{eqnarray*}
 \mathbb{K}[x_{1},\dots,x_{r}] & \leftrightarrow &\N^{r}\\
 x^{\mathbf{a}}=x_{1}^{a_{1}}\cdots x_{r}^{a_{r}} &\leftrightarrow & \mathbf{a}=(a_{1},\dots, a_{r}).
\end{eqnarray*}
Observe that $x^{\mathbf{a}}|x^{\mathbf{b}}\Leftrightarrow \mathbf{a}\leq \mathbf{b}$, where $\le$ denotes the coordinate-wise partial order in $\N^{r}$.

  We recall the following standard definitions from, for example, \cite{Cox-Little-Oshea-1992}.
 
 \begin{definition}
 A \emph{monomial ideal}  is an ideal $I \subset \mathbb{K}[x_{1},\dots,x_{r}]$ that can be generated by monomials. That is to say, $I$ consists of all polynomials which are finite sums of the form $\sum_{\alpha \in A}h_{\alpha}x^{\alpha}$, where $A\subset \N^{r}$ is a fixed finite subset of monomials and $h_{\alpha} \in \mathbb{K}[x_{1},\dots,x_{r}]$, for each $\alpha \in A$.
We write $I=\big<x^{\alpha}: \alpha \in A\big>.$
 \end{definition}

A monomial ideal is also a vector space with basis the set $M$ of monomials (i.e., elements of $\N^{r}$) that it contains. For this reason we sometimes call $M$ itself a monomial ideal.
The property that an $M\subset \N^{r}$ needs to have in order to be an ideal in this sense is that $v \in M$ implies $v+w \in M$ for every $w\in \N^r$. Equivalently, $M \subset \N^{r}$ is a monomial ideal if its complement $S=\N^r\setminus M$ (called \emph{its set of standard monomials} of the ideal $I$) satisfies 
\[
v+w \in S \Rightarrow v,w \in S.
\]
 Observe that this is equivalent to what condition (2) in the definition of minimum distance diagram asks for the image of the map $D$. That is:
 
 \begin{lemma}
 \label{lemma:monomial_MDD}
 A map $D:\mathbb{Z}_N\to \N^r$ is an MDD for $C_N(s_1,\dots,s_r)$ if, and only if, $D$ satisfies condition (1) in the definition and its image is the complement of a monomial ideal $M\subset \N^r$. 
 \end{lemma}

\begin{definition}
 A \emph{monomial ordering} on $\mathbb{K}[x_{1},\dots, x_{r}]$ is any relation $\prec$ on $\N^{r}$, or equivalently, any relation on the set of monomials $x^{\alpha}$, $\alpha \in \N^{r}$, satisfying:
 \begin{itemize}
 \item $\prec$ is a total ordering on $\N^{r}$.
 \item If $\alpha\prec \beta$ and $\gamma \in \N^{r}$, then $\alpha +\gamma \prec \beta+\gamma$ (in particular, but not only, $\prec$ extends the partial coordinate-wise order $\le$.
 \item $\prec$ is a well-ordering on $\N^{r}$. This means that every nonempty subset of $\N^{r}$ has a smallest element under $\prec$.
 \end{itemize}
 
A monomial ordering is \emph{graded} if it extends the (partial) ordering given by the $L_1$ norm (or total degree) of monomials.
 \end{definition}

 \begin{remark}
 \label{rem:coherent-MDD}
\rm
A graded monomial ordering $\prec$ in $\N^r$ induces an MDD of every $r$-tuple-loop network $C_N(s_1,\dots,s_r)$. Namely, the map:
\begin{eqnarray*}
D_\prec:\mathbb{Z}_{N} &\longrightarrow &\mathbb{N}^{r}\\
c &\longmapsto &\min{}_{\prec}\{(a_1,\dots,a_r)\in \N^{r}: c= a_{1}s_{1}+\cdots+a_{r}s_{r}\}.
\end{eqnarray*}
Indeed, the fact that $\prec$ is graded implies that the $\min_\prec$ in the formula is one of the solutions with minimum total degree (that is, with minimum $L_1$-norm). This implies that $D_\prec$ satisfies condition (1) of the definition of MDD. Condition (2) follows from Lemma~\ref{lemma:monomial_MDD}.
\end{remark}

 For practical purposes, one normally needs to know the ordering $\prec$ for monomials with $L_1$-norm bounded by a constant. For example, in order to construct $D_\prec$ in the above remark we will never need to compare monomials of $L_1$-norm bigger than $N$ (those will never give a minimum). An easy way of specifying such a ``bounded'' monomial ordering is via a ``weight vector'' $w=(w_1,\dots,w_r)\in [0,\infty)^r$. The ordering $\prec_w$ represented by it is  
  \[
{ \mathbf a} \prec_w {\mathbf b} \quad\Leftrightarrow\quad
{ \mathbf a} \cdot w  < {\mathbf b} \cdot w. 
 \]
Of course, $w$ has to be chosen ``sufficiently generic'' so that equality never arises in the equation for the 
finite (since they have bounded $L_1$-norm) set of vectors we are interested in . For example, the lexicographic ordering on $\mathbb K[x_1,\dots,x_r]$ is the ordering $\prec_w$ obtained when $w_i \gg w_{i+1}$ for every $i$.

The same applies if we want a graded monomial ordering. In this case, we define the ordering by first looking at the $L_1$-norm of the vectors ${ \mathbf a} $ and $ {\mathbf b}$, and using $w$ only to ``break ties''. The MDD of Remark~\ref{rem:coherent-MDD} for a graded order looked in this fashion has the following interpretation: the edges of type $s_i$ in the multi-loop network $C_N(s_1,\dots,s_r)$ have been assigned a weight $w_i$. The MDD chooses, among all the routes of minimal length, the one that has minimum weight.

\quad

Given an ideal $I \subset \mathbb{K}[x_{1},\dots,x_{r}]$ (not necessarily monomial) and a monomial ordering $\prec$ (not necessarily graded), it is well known that the leading monomials of the polynomials in $I$ form a monomial ideal, called the \emph{initial ideal of $I$ with respect to the ordering $\prec$}. The calculation of an initial ideal of $I$ is equivalent to that of a Gr\"obner basis.

The main result of \cite{Gutierrez-etal-2006} for our purposes is that the MDD's obtained by monomial orderings in Remark~\ref{rem:coherent-MDD} are, in fact, initial ideals of a certain ideal $I$ associated to $C_N(s_1,\dots,s_r)$.
\medskip

Let us recall that an \emph{integer lattice} is an additive subgroup of $\Z^{r}$ and that to every integer lattice $\mathcal{L} \subset \Z^{r}$ one can naturally associate the following \emph{lattice ideal}:
\[
I_{\mathcal{L}}:=\big < x^{\mathbf{a}^{+}}-x^{\mathbf{a}^{-}}: \mathbf{a} \in \mathcal{L} \big > \subset \mathbb{K}[x_{1},\dots,x_{r}]
\]
where $\mathbf{a}=\mathbf{a}^{+}-\mathbf{a}^{-}$ is the unique descomposition of $\mathbf{a}$ with $\mathbf{a}^{+}, \mathbf{a}^{-} \in \N^{r}$. A lattice ideal is a \emph{binomial ideal} (it is generated by binomials). Moreover (see \cite{Sturmfels-etal-1995}):
\[
x^{\mathbf{a}}-x^{\mathbf{b}} \in I_{\mathcal{L}} \Leftrightarrow \mathbf{a}-\mathbf{b} \in \mathcal{L}.
\]
\begin{theorem}[G\'omez et al., 2006 \cite{Gutierrez-etal-2006}]
\label{thm:Gutierrez-etal-2006}
Let $N, s_{1}, \dots, s_{r} \in \N$. Let us consider the lattice 
\[
\mathcal{L}:=\{(a_{1}, \dots, a_{r}) \in \Z^{r}: a_{1}s_{1}+\cdots+a_{r}s_{r}\equiv 0 \mod N\}.
\] 

Then, for every graded monomial ordering $\prec$, 
(the image of) the MDD $D_\prec$ of Remark~\ref{rem:coherent-MDD} coincides with (the set of standard monomials of) the initial ideal $In_\prec(I)$, where $I$ is the lattice ideal of 
$\mathcal{L}$.
\end{theorem}

The binomial ideal $I$ of the lattice $\mathcal{L}$
 in the theorem will be called the \emph{binomial ideal associated to the network $C_{N}(s_{1}, \dots, s_{r})$}.
G\'omez et al.~(\cite[Prop.~5]{Gutierrez-etal-2006}) show that $I$ can be generated by $r+2$ binomials, the calculation of which amounts to find integers $\lambda_1,\dots,\lambda_r,\mu$ satisfying
\[
\gcd(s_1,\dots,s_r, N) = \lambda_1 s_1 + \cdots + \lambda_r s_r + \mu N.
\]

Here we offer an alternative expression of $I$. We do not claim it to be algorithmically better (it is based in computing an eliminitaion ideal instead of a gcd) but it is theoretically ``more compact'' and easier to type-in in a computer algebra system, which is good for small examples where computation time is not an issue:

\begin{theorem}
\label{thm:elimination-ideals}
The binomial ideal associated to the network $C_{N}(s_{1}, \dots, s_{r})$ equals the elimination ideal of the variable $t$ in the following binomial ideal:
\[
\widetilde{I}:=\big <t^{N}-1,t^{s_{1}}-x_{1}, \dots, t^{s_{r}}-x_{r} \big >.
\]
\end{theorem}

In particular, if $\mathcal{G}$ is a reduced Gr\"obner basis, with respect to the elimination ordering, of the ideal $\widetilde{I}$, then the set of leading monomials of the elements of $\mathcal{G}\cap\mathbb{K}[x_1,\dots,x_r]$ constitutes a minimal system of generators of (the complement of) an MDD since (as we said before) the calculation of a Gr\"obner basis is equivalent to the calculation of an initial ideal.

Note that the ideal $\widetilde{I}$ is also the ideal of a lattice, namely the following one:
\[
\widetilde{\mathcal{L}}:=\big < Ne_{t},s_{1}e_{t}-e_{1}, \dots, s_{r}e_{t}-e_{r} \big >.
\]
This lattice is very special in the sense that the generators of the lattice directly give a system of generators of the ideal. In general, a lattice ideal may need more generators than the lattice.

\medskip

We finish this section with the observation that the reciprocal of Theorem~\ref{thm:Gutierrez-etal-2006} is not true. That is, there are multi-loop networks with MDD's that cannot be derived from monomial orderings as in Remark~\ref{rem:coherent-MDD} or Theorem~\ref{thm:Gutierrez-etal-2006}.

\begin{example}
\label{ex:non-coherent}
In the network $C_{8}(1,3,5,7)$ the set of monomials $M=\{x_{1}, x_{2},$ $x_{3},x_{4},x_{1}^{2},x_{2}^{2},x_{3}x_{4}\}$ ``is" an MDD but it is not $D_\prec$ for any ordering.

\rm
Indeed, the monomials $x_{1}, x_{2},$ $x_{3},x_{4}$ are in any MDD of $C_{8}(1,3,5,7)$ because the unique shortest path to the vertices $1$, $3$, $5$ and $7$ is via a single edge. For each of the other three vertices $2$, $4$ and $6$ there are two or three minimum paths, all of length two:
\begin{itemize}
\item $x_{1}^{2}$ and $x_{3}^{2}$ (and also $x_{2}x_{4}$, but we do not use it) for vertex $2$.
\item $x_{2}^{2}$ and $x_{4}^{2}$ (and also $x_{1}x_{3}$, but we do not use it) for vertex $6$.
\item  $x_{3}x_{4}$ and $x_{1}x_{2}$ for vertex $4$.
\end{itemize}
The  $3 \times 3 \times 2=18$ possible choices of minimum paths give 18 different MDD's. In $M$ we have chosen the first path described for each vertex. This choice is not compatible with any monomial ordering $\prec$ because:
\begin{itemize} 
\item If $\prec$ selects $x_{1}^{2}$ from the binomial $x_{1}^{2}-x_{3}^{2}\in I$ then $x_1 \prec x_3$.
\item If $\prec$ selects $x_{2}^{2}$ from the binomial $x_{2}^{2}-x_{4}^{2}\in I$ then $x_2 \prec x_4$.
\item Hence, $x_1x_2 \prec x_3x_4$ and the monomial selected in $x_{3}x_{4} - x_{1}x_{2}$ should have been $x_1x_2$.
\end{itemize}
\end{example}

Mimicking the literature on $A$-graded ideals and toric Hilbert schemes (see Chapter 10 in \cite{Sturmfels-1996}) we call MDD's \emph{coherent} or \emph{non-coherent} depending on whether they can be obtained from monomial orderings  or not. 

If the reader goes back to the proof of Lemma~\ref{lemma:2MDDs} he or she will notice that it is based on the facts that there are only two graded monomial orderings in two variables, and 
 double loop networks do not have non-coherent MDD's. 

For triple-loop networks we do not know whether non-coherent MDD's exist. But, even if they do not, in the next section we show networks with arbitrarily many coherent MDD's. A crucial object in our construction, implicit also in Example~\ref{ex:non-coherent} and in the proof of Lemma~\ref{lemma:2MDDs}, is the following \emph{homogeneous sublattice} $\mathcal{L}_0$ of the lattice $\mathcal{L}$ of $C_N(s_1,\dots,s_r)$:
\[
\mathcal{L}_{0} := \mathcal{L}\cap \{(x,y,z)\in \mathbb{Z}^{3}:\, x+y+z=0\}.
\]
It is clear that the lattice $\mathcal{L}_{0}$ is the source of non-uniqueness of MDD's, as is explicit in the following result:

\begin{lemma}[\cite{Gutierrez-etal-2006}]
\label{lemma:non-uniqueness}
Let $M\subset \N^r$ be an MDD for the network $C_N(s_1,\dots,s_r)$. Then, $M$ is the unique MDD for that network if and only if there is no $\mathbf a\in M$ and $\mathbf b\in \mathcal{L}_0\setminus\{0\}$ such that $\mathbf a+\mathbf b\in \N^r$.
\end{lemma}

\begin{proof}
Observe that, for each $\mathbf a\in M$, the paths in $C_N(s_1,\dots,s_r)$ that lead to the same vertex as $\mathbf a$ correspond precisely to the vectors in $(\mathbf a+\mathcal{L}_0) \cap \N^r$. In particular, if that intersection contains only $\mathbf a$ for each $\mathbf a\in M$, then $M$ is the unique MDD.

Reciprocally, suppose that for some $\mathbf a\in M$ and for some non-zero $b\in \mathcal{L_0}$ we have 
$\mathbf a':=\mathbf a+ \mathbf b\in \N^r$. Then, consider any graded monomial ordering $\prec$ such that $\mathbf a'\prec \mathbf a$ (for example, a degree-lexicographic ordering starting with any variable whose coordinate is bigger in $\mathbf a$ than in $\mathbf a'$). Then, $x^{\mathbf a}$ is in the initial ideal $\in_\prec(I)$ since it is the leading monomial of
$x^{\mathbf a} - x^{\mathbf a'}\in I$. Hence, $\mathbf a$ is not in the coherent MDD produced by $\prec$.
\end{proof}

That is, 
choices in the construction of an MDD correspond to elements of $\mathcal{L}_{0}$.

 \section{Coherent MDD's and Hilbert bases of lattice cones}
 \label{sec:ManyMDD's}

As we said after Remark~\ref{rem:coherent-MDD}, every coherent minimum distance diagram for  a given network $C_N(s_1,\dots,s_r)$ 
(more generally, every initial ideal for a given ideal $I\subset {\mathbb K}[x_1,\dots,x_r]$ is the MDD constructed
from a sufficiently generic weight vector $w=(w_1,w_2,\dots,w_r)\in \R^r$ in the following fashion:
$D_w(i)$ is \emph{the} ${\mathbf a}\in \N^{r}$ that minimizes $w\cdot {\mathbf a}$ among those that 
minimize $||{\textbf a}||$. 

 It is clear that the definition of $D_w$ is not affected  when we multiply $w$ by a positive constant, or when we add to it a real multiple of $(1,\dots,1)$. Hence, \emph{coherent MDD's are parametrized by rays in the hyperplane }
 \[
 H_0=\{(w_1,\dots,w_r)\in \R^r: w_1+\cdots+w_r=0\}.
 \]
 
From now on we assume that $r=3$, so that $H_0$ is a 2-dimensional plane. 
If a sufficiently generic $w$ produces a certain MDD $D_w$ and we perturb it to a very very close $w'$, the new monomial ordering $\prec_{w'}$ will be the same as $\prec_w$ (not over all $\N^3$ but over the bounded, finite, part of $\N^3$ that is of interest once $C_N(s_1,\dots,s_3)$ has been fixed). The regions of $H_0$ corresponding to vectors $w$ that produce the same graded order are two-dimensional open and rational cones, each bounded by two rays, and cyclically ordered around the origin in $H_0$. That is, they form a 2-dimensional complete polyhedral fan. 
It is obvious that the number of MDD's (regions in the fan) equals the number of rays between 
consecutive regions. Our goal in this section is to characterize those rays.

So, for the rest of this section, let $w_0\in H_0$ be a non-zero vector in the common boundary of two cones, and let $w_+$ and $w_-$ be two sufficiently small perturbations of it, each lying in the interior of one of the adjacent cones. Let $D_+=D_{w_+}$ and $D_-=D_{w_-}$ be the MDD's produced by $w_+$ and $w_-$, respectively. 

Crucial in our characterization is going to be the homogenous lattice $\mathcal{L}_0$ of $C_N(s_1,s_2,s_3)$, and at the \emph{Hilbert bases} of its intersection with orthants.
Recall that $\mathcal{L}_0$ is:
\[
\mathcal{L}_{0} := \{(a_1,a_2,a_3)\in \mathbb{Z}^{3}:\, a_1 s_1 + a_2 s_2 + a_3 s_3 =0 \, (\mod N), \text{ and } a_1 + a_2 +a_3=0\}.
\]

\begin{lemma}
There is a non-zero $\mathbf a \in \mathcal{L}_0$ such that $w_0\cdot \mathbf{a} =0$.
\end{lemma}

\begin{proof}
Let $i\in \Z_N$ be such that $D_+(i)\ne D_-(i)$. Let $ \mathbf{a}= D_+(i) - D_-(i)$.
Then, $D_+(i)$ and $D_-(i)$ represent routings to the same vertex $i$ of $C_N(s_1,s_2,s_3)$, and of the same length, so $ \mathbf{a} \in \mathcal L_{0}$. 
On the other hand, $D_+(i) \prec_{w_+} D_-(i)$ and 
$D_-(i) \prec_{w_-} D_+(i)$, that is, 
\[
w_+\cdot \mathbf{a} <0,
\qquad
\text{and}
\qquad
w_-\cdot \mathbf{a} >0.
\]
By continuity, $w_0\cdot \mathbf{a} =0$.
\end{proof}

For each choice of signs $\epsilon =(\epsilon_1, \epsilon_2, \epsilon_3)\in \{-,+\}^3$ we consider the semigroup 
orthant
\[
S^\epsilon:= \mathcal{L}_0 \cap \{(a_1,a_2,a_3)\in \mathbb{Z}^{3}:\, \epsilon_i a_i \ge 0, \forall i\}.
\]
The \emph{Hilbert basis} of $S^\epsilon$ is the set of elements of $S^\epsilon$ that cannot be expressed as a sum of two non-zero elements of it. That is, it is the unique minimal generating system of  $S^\epsilon$ as a semigroup.

Let $\mathbf a \in \mathcal L_0\setminus \{(0,0,0)\}$  have minimum norm among the vectors satisfying $w_{0}\cdot \mathbf a=0$. In part (3) of the following result it is crucial to assume that $\mathbf a$ has only one negative entry. 
This is no loss of generality since it can be achieved  by changing $\mathbf a$ to $-\mathbf a$, if necessary. But the same would not be true for $r>3$, so only the first two parts in the lemma generalize to arbitrary $r$.

\begin{lemma}
\label{lemma:hilbert}
Let $ \mathbf{a}= \mathbf{a}_+  - \mathbf{a}_-$ be the unique decomposition of $ \mathbf{a}$ into two non-negative vectors, where $ \mathbf{a}$ has minimal norm among  the elements of $\mathcal L_0$ with  $w_0\cdot \mathbf{a} =0$. Then,
\begin{enumerate}
\item $\mathbf{a}_+$ and $\mathbf{a}_-$ lie, respectively, in the two MDD's $D_+ $ and $ D_-$ ``incident to''
$w_0$. In particular, they
 represent minimum routings in $C_N(s_1,s_2,s_3)$. 

\item ${\mathbf a}$ is in the Hilbert basis of $S^{\epsilon}$, where $S^{\epsilon}$ is any orthant semigroup containing ${\mathbf a}$.

\item If $\mathbf a$ has only one negative entry then, 
for every $ \mathbf{b}\in S^{\epsilon}$ with 
$|| \mathbf{b}||\le || \mathbf{a}||$ we have $w_0\cdot \mathbf{b} \ge0$.
\end{enumerate}
\end{lemma}

\begin{proof}
As in the previous lemma, let  $i\in \Z_N$ be such that $D_+(i) \ne D_-(i)$, which implies that $w_0\cdot (D_+(i) - D_-(i)) = 0$. Hence, $D_+(i) -  D_-(i)$ is proportional to $ \mathbf{a}$. By exchanging $D_+$ and $D_-$ if necessary, there is no loss of generality in assuming that $D_+(i) - D_-(i)$ is a positive multiple of $ \mathbf{a} =\mathbf{a}_+  - \mathbf{a}_- $, which implies that $D_+(i)$ and $D_-(i)$ are respective positive multiples of $\mathbf{a}_+$ and  $ \mathbf{a}_- $, with the same factor. This factor must be an integer, by minimality of $\mathbf a$, and then the second condition in the definition of an MDD implies that $\mathbf{a}_+$ and  $ \mathbf{a}_- $ represent also routings in  $D_+ $ and $ D_-$, respectively.

For part (2), suppose that ${\mathbf a}$ was not in the Hilbert basis. That is,
let ${\mathbf a}= {\mathbf b} + {\mathbf c}$, with  ${\mathbf b}, {\mathbf c}\in S^{\epsilon}\setminus\{(0,0,0)\}$.
Let $ \mathbf{b}= \mathbf{b}_+  - \mathbf{b}_-$ and $ \mathbf{c}= \mathbf{c}_+  - \mathbf{c}_-$ be the decompositions of $ \mathbf{b}$ and $ \mathbf{c}$ into positive and negative parts. Observe that
\[
 \mathbf{a}_+= \mathbf{b}_+  + \mathbf{c}_+= \mathbf{b}  + (\mathbf{c}_+ + \mathbf{b}_-)= 
 (\mathbf{b}_+  + \mathbf{c}_-) + \mathbf{c}.
 \]
In particular, as paths in the network, $\mathbf{b}_+  + \mathbf{c}_-$ and $\mathbf{c}_+ + \mathbf{b}_-$ lead to the same vertex as $\mathbf{a}_+$. By part (1), then,
\begin{eqnarray*}
w_0 \cdot  (\mathbf{b}_+  + \mathbf{c}_-) &\ge& w_0 \cdot \mathbf{a}_+ = w_0 \cdot \mathbf{a}_-,\\
w_0 \cdot  (\mathbf{b}_-  + \mathbf{c}_+) &\ge& w_0 \cdot \mathbf{a}_+ = w_0 \cdot \mathbf{a}_-.
\end{eqnarray*}
This, together with the previous equalities implies
\[
w_0\cdot \mathbf{b}=0, \qquad  w_0\cdot \mathbf{c} = 0,
\]
which contradicts the minimality in the choice of ${\mathbf a}$.

For part (3), to fix notation assume, without loss of generality, that $\epsilon= (-,+,+)$. We can then write
$ \mathbf{a}=(-a_1,a_2,a_3)$ with $a_1,a_2,a_3 \ge 0$.
Let $ \mathbf{b}=(-b_1,b_2,b_3)\in S^{(-,+,+)}$ have 
$|| \mathbf{b}||\le || \mathbf{a}||$. In particular, $b_1=|| \mathbf{b}||/2 \le a_1$, so that 
$\mathbf{a}_- + \mathbf{b} \in \N^3$. Since $\mathbf{b}\in \mathcal L$ and since $D_-(i)=\mathbf{a}_-$, we have
\[
\mathbf{a}_- \prec_{w_-} \mathbf{a}_- + \mathbf{b} 
\qquad
\Rightarrow
\qquad
w_-\cdot \mathbf{b} >0
\qquad
\Rightarrow
\qquad
w_0\cdot \mathbf{b} \ge0.
\]
The latter implication is by continuity.
\end{proof}

Perhaps more interestingly, we also have the following converse to this lemma:

\begin{lemma}
\label{lemma:hilbert2}
Let $w_0\in H_0\setminus \{0\}$ and let $\mathbf a= \mathbf{a}_+  - \mathbf{a}_-$ satisfy all the hypotheses of Lemma~\ref{lemma:hilbert}. That is to say:
\begin{enumerate}
\item[(0)] $\mathbf a$  has minimum norm among the elements of $\mathcal L_0$ orthogonal to $w_0 $.
\item $\mathbf{a}_+$ and $\mathbf{a}_-$ represent minimum routings in $C_N(s_1,s_2,s_3)$. 
\item ${\mathbf a}$ is in the Hilbert basis of the corresponding $S^{\epsilon}$.
\item $\epsilon$ has a single negative entry and for every $ \mathbf{b}\in S^{\epsilon}$ with 
$|| \mathbf{b}||\le || \mathbf{a}||$ we have $w_0\cdot \mathbf{b} \ge0$.
\end{enumerate}

Then, $w_0$ is the common boundary ray of two MDD cones in $H_0$.
\end{lemma}

\begin{proof}
The only thing we need to prove is that $w_0 \cdot \mathbf c > w_0 \cdot \mathbf a_-$ for every 
$\mathbf c \in \mathbb N^3$ different from $\mathbf a_+$ and $\mathbf a_-$ and with the same norm, and leading to the same vertex $i$ of the network. Indeed, if this is the case, every sufficiently small perturbation $w'$ of $w_0$ will select either 
$\mathbf a_+$ or $\mathbf a_-$ as the path to choose for the MDD. Which one is selected will only depend on the sign of $w'\cdot \mathbf a$.

So, let $\mathbf c$ be in that conditions. Observe that, then, $||\mathbf c||=||\mathbf a_-||$ (by part (1)) and hence
$\mathbf c - \mathbf a_- \in \mathcal L_0$.
The fact that $\mathbf a_-$ has a unique non-zero entry implies that 
$\mathbf c - \mathbf a_-$ is in $S^\epsilon$ and that it has the same or smaller norm as $\mathbf a$. By part (2), then,
$w_0\cdot (\mathbf c - \mathbf a_-)\ge 0$. Equality is impossible, since it would imply that 
 $\mathbf c - \mathbf a_-$ is proportional to $\mathbf a$, in violation with the minimality of $||\mathbf a||$.
 Hence, $w_0\cdot (\mathbf c - \mathbf a_-) > 0$ and  $w_0\cdot \mathbf c > w_0 \mathbf a_-$, as we wanted to proof.
\end{proof}

Observe that, actually, in this proof we do not use that $\mathbf a$ is in the Hilbert basis. But, as we saw in the proof of Lemma~\ref{lemma:hilbert}, that property follows from (0) and (1).

Lemma~\ref{lemma:hilbert} can be read in reverse: by part (2), every ray incident to two cones of the fan of MDD's is orthogonal to an element $\mathbf a$ in the Hilbert basis of one of the semigroups  $S^{\epsilon}$. Of course, $S^{-\epsilon}= -S^\epsilon$. Since, also, $S^{(+,+,+)}=S^{(-,-,-)}=\{0\}$, there is no loss of generality in considering only the three 
semigroup orthants $S^{(-,+,+)}$, $S^{(+,-,+)}$ and $S^{(+,+,-)}$ to which part (3) applies. With this we get:

\begin{corollary}
\label{cor:hilbert}
The number of coherent MDD's for a network with homogeneous lattice $\mathcal L_0$ is bounded above by the sum of cardinalities of the Hilbert bases of the three octant semigroups $S^{(-,+,+)}$, $S^{(+,-,+)}$ or $S^{(+,+,-)}$.
\end{corollary}

\begin{proof}
Each element $\mathbf a$ in one of the three Hilbert bases can in principle produce two rays $w$ (the two rays orthogonal to  $\mathbf a$. This in principle allows for twice the number of MDD's that we want to prove. But:
\begin{itemize}
\item If $||\mathbf a||$ is not the minimum among the norms of non-zero elements of its semigroup $S^\epsilon$, then only one of the two rays orthogonal to $||\mathbf a||$ satisfies part (3) of Lemma~\ref{lemma:hilbert}.
\item If $S^\epsilon$ has several non-zero elements $\mathbf a_1, \mathbf a_2, \dots, \mathbf a_k$ with minimum norm, then in total there are two rays orthogonal to one of them and satisfying condition (3): the interior normals of the cone $\pos(\mathbf a_1, \mathbf a_2, \dots, \mathbf a_k)$.
\item Only if $\mathbf a$ is the unique element with minimum norm among non-zero elements of its semigroup $S^\epsilon$, then the two rays orthogonal to $||\mathbf a||$ satisfy part (3).
\end{itemize}

Thus, only for three of the elements in the union of the Hilbert bases we can get two rays. But there are also  three rays that are counted twice in this process. Indeed,  the ray generated by $(-1,-1,2)$ arises both from the Hilbert basis element $(a,-a,0)\in S^{(+,-,+)}$ and from its opposite $(-a,a,0)\in S^{(-,+,+)}$, and the same happens for the rays generated by $(2,-1,-1)$ and $(-1,2,-1)$.
\end{proof}

More interesting than the statement of this corollary is the explicit way described in its proof to get a list of rays susceptible of being incident to two MDD's. Let us see this in two examples. The second one also shows that 
the bound in this corollary is not tight for every network. The reason is that this bound takes only $\mathcal L_0$ into account, while the fan of MDD's does not only depend on $\mathcal L_0$ (as is implicit also in part (1) of Lemma~\ref{lemma:hilbert}).

\begin{example}
\rm
Consider the lattice $\mathcal L_0$ of Figure~\ref{fig:9MDD's}, generated by (for example), the vectors $(3,0,-3)$ and $(1,1,-2)$. In this and the following pictures the blue dots represent the elements of $\mathcal L_0$, and the white dots the rest of integer points in the plane $x+y+z=0$. Only the parts in the three octants that we need to study are shown, and the black dots represent the Hilbert basis of each. The following is the list of the nine Hilbert basis elements and the rays orthogonal to them that satisfy condition (2) of the lemma. As predicted in the proof of Corollary~\ref{cor:hilbert}, three of them arise twice in the list:
\begin{figure}[htb]
     \centering
      \includegraphics[width=6 cm]{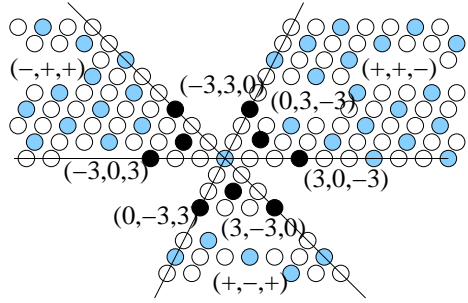}
    \caption{Lattice $\mathcal{L}_{0}$ of $C_{9}(1,4,7)$}
    \label{fig:9MDD's} 
  \end{figure}
\[
\begin{matrix}
\mathbf a \in S^\epsilon & & w_0\\
(0,3,-3) \in S^{(+,+,-)} & \to  &(2,-1,-1) \\
(0,-3,3) \in S^{(+,-,+)}  & \to  &(2,-1,-1) \\
(3,0,-3) \in S^{(+,+,-)}  & \to  &(-1,2,-1) \\
(-3,0,3)  \in S^{(-,+,+)} & \to  &(-1,2,-1) \\
(3,-3,0)  \in S^{(+,-,+)} & \to  &(-1,-1,2) \\
(-3,3,0)  \in S^{(-,+,+)} & \to  &(-1,-1,2) \\
(-2,1,1)  \in S^{(-,+,+)} & \to  &(0,1,-1) ,\ (0,-1,1)  \\
(1,-2,1)  \in S^{(+,-,+)} & \to  &(1,0,-1) ,\ (-1,0,1)  \\
(1,1,-2)  \in S^{(+,+,-)} & \to  &(1,-1,0) ,\ (-1,1,0)  \\
\end{matrix}
\]
Hence, for every network $C_N(s_1,s_2,s_3)$ having this lattice we have at most nine coherent MDD's. The bound is tight since it is achieved for the network $C_9(1,4,7)$, as we saw in Figure~\ref{fig:3Pasosy9MDDs}. But other networks with the same homogeneous lattice may have strictly less coherent MDD's. For example, the network $C_6(1,3,5)$ has only four (coherent or not) MDD's: there are two choices of path to vertex 2, and two choices to vertex 4.
\end{example}

\begin{example}
\label{exm:2min}
\rm
Consider now the lattice $\mathcal L_0$ of Figure~\ref{fig:2min}. In the semigroup $S^{(-,+,+)}$ there are two Hilbert basis elements with minimal norm, namely $(-4,3,1)$ and $(-4,1,3)$. Hence, the count of Corollary~\ref{cor:hilbert} still has an excess of one: there are 10 Hilbert basis elements in total but only 9 rays susceptible of being incident to two MDD's. 

If we think of this lattice as the $\mathcal L_0$ of the network $C_8(2,3,7)$ we easily see that the number of MDD's is merely two: the minimum paths in the network to the vertices $2$, $3$, $7$, $1=7+2$, $4=2+2$ and $5=3+2$ are unique. Our only choice is in the minimum path to $6=3+3=7+7$.

However, this $\mathcal L_0$ is also the homogeneous lattice of the network $C_{72}(19,28,64)$, and in this one we do get the 9 MDD's allowed by Lemma~\ref{lemma:hilbert}. That this is the general situation is proved in Theorem~\ref{thm:N'} below; see in particular,  Example~\ref{exm:2min-cont}.
\begin{figure}[htb]
     \centering
      \includegraphics[width=9cm]{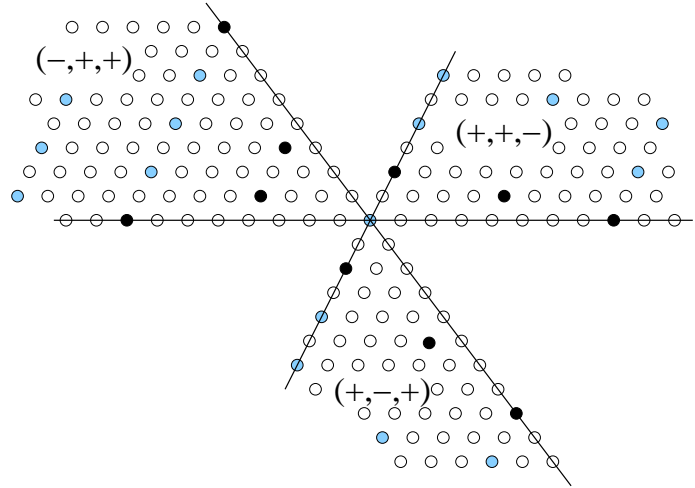}
    \caption{Lattice $\mathcal{L}_{0}$ of $C_{8}(2,3,7)$}
    \label{fig:2min} 
  \end{figure}
\end{example}

\begin{remark}
\rm
To better understand the examples, observe that for a semigroup $S=\mathcal L \cap C\setminus \{(0,0)\}$ obtained as the intersection of a 2-dimensional lattice with a linear cone, the Hilbert basis of $S$ coincides with the elements in the boundary of its lower hull. That is:
\[
\operatorname{Hilb}(S)=\{\mathbf a\in S : \forall \lambda <1, \quad \lambda \mathbf a \not\in \operatorname{conv}(S)\}.
\]

Indeed, 
if $\mathbf a = \mathbf b + \mathbf c$ is not a Hilbert basis element, then $\mathbf a/2$ is the midpoint of the segment $\mathbf b  \mathbf c$, so that $\mathbf a$ is not in the lower hull of $S$ (this implication holds in every dimension). 

Conversely, suppose that $\mathbf a$ is not in the lower hull of $S$. Let $\mathbf b$ and $\mathbf c$ be consecutive elements of $S$ in its lower hull and such that $\mathbf a \in \pos (\mathbf b, \mathbf c)$. That is, 
$\mathbf a = \lambda \mathbf b + \mu \mathbf c$ for nonnegative real numbers $\lambda$ and $\mu$.
Then, by construction, the triangle $O\mathbf b \mathbf c$ contains no points of $\mathcal L$ other than its vertices. This (for example by Pick's Theorem) implies that $\mathbf b$ and $\mathbf c$ are a lattice basis of $\mathcal L$, so that $\lambda$ and $\mu$ are integers and $\mathbf a$ is not in the Hilbert basis.
\end{remark}

\section{Triple loop networks with many MDD's}

In the previous section we have proved an upper bound of the number of coherent MDD's in terms of the homogeneous lattice $\mathcal L_0$. The goal of this section is two-fold:
\begin{enumerate}
\item Construct lattices where the bound is arbitrarily big (and which are lattices of some triple-loop network).
\item Show that the bound is attained: For every such lattice  there is some triple loop network with that homogeneous lattice and with that many coherent MDD's.
\end{enumerate}

We start with the second goal:

\begin{theorem}
\label{thm:N'}
Let $\mathcal L_0$ be the homogeneous lattice of some triple loop network $C_N(s_1,s_2,s_3)$.
Then, there is another triple loop network $C_{N'}(s'_1,s'_2,s'_3)$ which has the same homogeneous lattice and with the following property: if $w_0$ and $\mathbf a$ satisfy properties (0), (2) and (3) of Lemma~\ref{lemma:hilbert2} (which depend only on $\mathcal L_0$) then they also satisfy property (1).
\end{theorem}

Observe that not every sublattice  $\mathcal L$ (respectively, $\mathcal L_0$) of finite index in $\Z^r$ (respectively, in ${\mathcal Z}_0 = \{(a_1,\dots,a_r)\in \Z^r : \sum a_i=0\}$) is the lattice (respectively, the homogeneous lattice) of a multi-loop network. This happens if and only if the quotient groups $\Z^r / \mathcal L$ and $\mathcal Z_0 / \mathcal L_0$ are cyclic.

\begin{proof}
The proof has the following ingredients:
\begin{itemize}
\item For any $t,k\in \N$ such that $\gcd(k,N)=1$, 
 the following transformation on the triple loop network preserves the homogeneous lattice:
\[
N'=Nk,\qquad
s'_1 = t+ks_1, \qquad
s'_2 = t+ks_2, \qquad \text{and} \qquad
s'_3 = t+ks_3.
\]

Indeed, for a vector $(a_1,a_2,a_3)$ with $a_1+a_2+a_3=0$, the equation 
\[
a_1s'_1+a_2s'_2+a_3s'_3 = 0\ (\mod N')
\]
that defines the homogeneous lattice of $C_{N'}(s'_1,s'_2,s'_3)$ is equivalent to the equation 
\[
a_1 k s_1+a_2 k s_2+a_3 k s_3  = 0\ (\mod kN).
\]
The assumption that $k$ is prime with $N$ then allows us to remove the factor $k$ on both sides of this last equation.

\item If, moreover, $\gcd(k,t)=1$ and $k> ||\mathbf a_+||$, for a certain Hilbert basis element 
 $\mathbf a= \mathbf a_+ - \mathbf a_-$, then  $ \mathbf a_+$ and $ \mathbf a_-$ represent minimal routings in $C_{N'}(s'_1,s'_2,s'_3)$. 

Indeed, 
if we let $i\in \Z_{N'}$ be the vertex $(s'_1,s'_2,s'_3)\cdot \mathbf a_+ = (s'_1,s'_2,s'_3)\cdot \mathbf a_-$ 
to which these paths go
and we let $j$ be another vertex obtained by a shorter path of length, say, $l< ||\mathbf a_+||$,
we have that
$i = ||\mathbf a_+|| t (\mod k)$, while $ j = l t  (\mod k)$. Our assumptions imply that then $i\ne j$ because
$i =j$ and $\gcd(k,t)=1$ would imply $l= ||\mathbf a_+|| \mod k$, impossible since $l<||\mathbf a_+||< k$.

\end{itemize}
Thus, it suffices to let $k$ and $t$ be such that the assumptions in these properties hold for every $\mathbf a$ in the Hilbert basis. For example, it is easy to prove 
that $||\mathbf a_+||< N$ for every $\mathbf a$, so that taking $k=N+1$ and $t=1$ will do the job.
\end{proof}

\begin{example}[Example~\ref{exm:2min} continued]
\rm
\label{exm:2min-cont}
Let us look again at the lattice of Example~\ref{exm:2min}, in which there are nine rays $w$ satisfying conditions (0), (2) and (3) of Lemma~\ref{lemma:hilbert2}. As said there, this is the homogeneous lattice of the network $C_8(2,3,7)$, but this network has only two, instead of nine, MDD's. Applying to this network the procedure in the proof of Theorem 5.1, with $k=9$ and $t=1$ all the requirements in the proof are satisfied. Observe that $k=9=N+1$ is the minimum possible value that makes the proof work in this example, since $\mathbf a=(-8,0,8)$ is a Hilbert basis element in $S^{(-,+,+,)}$ with $||\mathbf a_-||=8$.
\end{example}

We finally show examples of lattices $\mathcal L_0$ with arbitrarily many Hilbert basis elements:

\begin{theorem}
\label{thm:many}
Let $q \in \mathbb{N}$, with $q-1$ not a multiple of three.  Let $N=1+q+q^2$
(so that $\gcd(q-1,N)=\gcd(q-1,3)=1$). Consider the the triple-loop network $C_{N}(1,q,q^{2})$. Then:
\begin{enumerate}
\item Its homogeneous lattice is symmetric under cyclic permutation of the three coordinates, and has $q+2$ Hilbert basis elements in each of the octants $S^{(-,+,+)}$, $S^{(+,-,+)}$, and $S^{(+,+,-)}$, namely (for the first one):
\[
(-N,0,N),  (-q-1,1,q)+ i (-q,q+1,-1),\qquad i=0,\dots,q. 
\]
\item For each of them there is a unique ray $w$ satisfying conditions (0), (2) and (3) of Lemma~\ref{lemma:hilbert2}.
\item As a consequence, the triple-loop network $C_{Nk}(t+k,t+qk,t+q^{2}k)$ has exactly $\mathbf{3(q+2)}$ coherent MDD's, for any $k$ bigger than $N$ and with $\gcd(t,k)=\gcd(k,N)=1$ (for example, $k=N+1$ and $t=1$).
\end{enumerate}
\end{theorem}

\begin{figure}[htb]
     \centering
      \includegraphics[width=10cm]{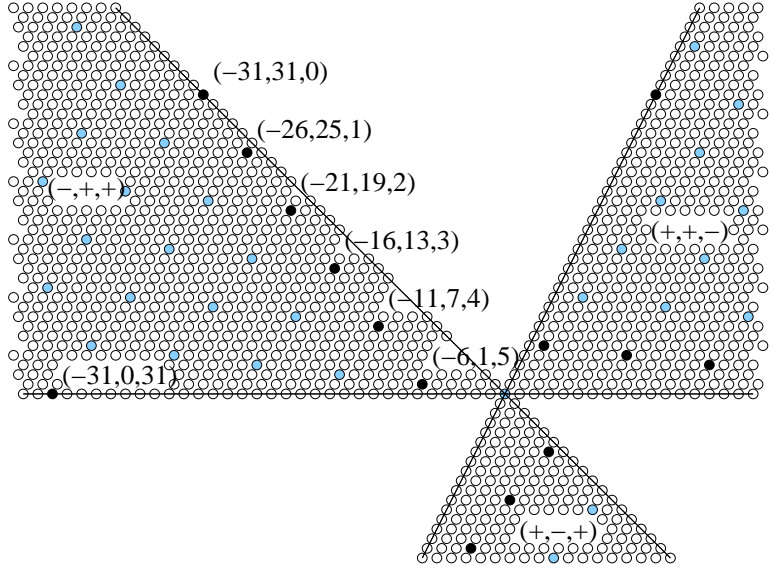}
    \caption{The lattice $\mathcal{L}_{0}$ of Theorem~\ref{thm:many}, with $q=5$}
    \label{fig:21MDD's} 
  \end{figure}

\begin{remark}
\rm
Observe that the condition $\gcd(t,k)=1$ is clearly necessary  for the network $C_{Nk}(t+k,t+qk,t+q^{2}k)$ to be connected. Together with $\gcd(q-1,N)=1$ it is also sufficient.
\end{remark}

\begin{proof}
 
Starting with the equation for the lattice  $\mathcal{L}$
\[
x+ q y+ q^{2} z\equiv 0 \, (\mod N),
\]
and using  $x+y+z=0$ to eliminate either one of the variables $x$, $y$ or $z$ we get the following three descriptions of 
the homogeneous lattice $\mathcal{L}_0$:
\[
\mathcal{L}_{0}=\left \{ {(q-1) y+(q^2-1) z\equiv 0 \, (\mod N) \atop x+y+z=0} \right. 
\]
\[
\mathcal{L}_{0}=\left \{ {(1-q) x+(q^2-q) z\equiv 0 \, (\mod N) \atop x+y+z=0} \right. 
\]
\[
\mathcal{L}_{0}=\left \{ {(1-q^2) x+(q-q^2) y\equiv 0 \, (\mod N) \atop x+y+z=0} \right. 
\]
Since $\gcd(q-1,N)=1$ we can divide by  $q-1$. This gives:  
\[
\mathcal{L}_{0}=\left \{ {y+(q+1)z\equiv 0 \, (\mod N) \atop x+y+z=0} \right. 
\]
\[
\mathcal{L}_{0}=\left \{ {qz - x\equiv 0 \, (\mod N) \atop x+y+z=0} \right. 
\]
\[
\mathcal{L}_{0}=\left \{ {(1+q) x +  q y\equiv 0 \, (\mod N) \atop x+y+z=0} \right. 
\]
Now, we divide the second and third equations by $q$ and $-(1+q)$ respectively, which can be done since 
$q^{-1}=-(q+1) \mod N$. This gives the following symmetric descriptions, which prove part (1) of the statement:
\[
\mathcal{L}_{0}=\left \{ {y+(q+1)z\equiv 0 \, (\mod N) \atop x+y+z=0} \right. 
\]
\[
\mathcal{L}_{0}=\left \{ {z + (q+1)x\equiv 0 \, (\mod N) \atop x+y+z=0} \right. 
\]
\[
\mathcal{L}_{0}=\left \{ {x +  (q+1) y\equiv 0 \,(\mod N) \atop x+y+z=0} \right. 
\]

For the rest of the proof we concentrate in the octant  $S^{(-,+,+)}$. We first prove that the $q+2$ vectors stated are in the Hilbert basis. For $(-N,0,N)$ this is obvious: any vector $(- a,0,a)$ in $\mathcal L_0$ will have $a\cdot (q+1)=0 \, (\mod N)$, that is, $a=0 \, (\mod N)$. For the rest we observe that for any element $(x,y,z) \in S^{(-,+,+)}\setminus{(0,0,0)}$ we have that $y + (1+q) z$ is positive, and a multiple of $N$. Hence, all those with $y + (1+q) z=N$ must be in the Hilbert basis. It is easy to check that those are precisely the vectors of the form
\[
(x,y,z)= (-q-1,1,q)+ i (-q,q+1,-1),\qquad i=0,\dots,q
\]

That there are no other elements in the Hilbert basis can be proved as follows: indeed, let $\mathbf b=(-b_2-b_3,b_2, b_3) \in S^{(-,+,+)}$ be such that
\[
b_2 + (1+q) b_3 \ge 2N.
\]
We distinguish three cases:
\begin{itemize}
\item If $b_2=0$, then the only possibility is $\mathbf b = (-N,0,N)$.
\item If $b_2>0$ and $b_3 \ge q$ we can write
\[
\mathbf b=(-q-1,1,q) + (-b_2-b_3+q+1,b_2-1, b_3-q),
\]
which proves that $\mathbf b$ is not in the Hilbert basis.
\item If $b_3 < q$, then $b_2 >N$ and we can write
\[
\mathbf b=(-N,N,0) + (-b_2-b_3 + N ,b_2-N, b_3),
\]
which also proves that $\mathbf b$ is not in the Hilbert basis.
\end{itemize}
This finishes the proof of part (2). Part (3) is a direct application of (the proof of) Theorem~\ref{thm:N'}.

\end{proof}

\subsection*{Acknowledgement:} 
Almost all of the figures of Sections 1 and 2 have been drawn automatically with a program  created by \'Alvar Ibeas. We thank him for this software which we used not only for the final pictures but also for the exploration and understanding of MDD's in different examples of networks.

\bibliographystyle{abbrv}

\end{document}